\documentclass[12pt]{article}
\usepackage[T2A]{fontenc}
\usepackage{amsmath}
\usepackage[cp1251]{inputenc}

\begin{document}
{\Large

\vskip 0.1cm

\begin{center}
ON THE RECONSTRUCTION OF THE MATRIX \\ FROM ITS MINORS \vskip 0.2cm
{\bf Mouftakhov~A.~V.(Ramat-Gan (Israel), BIU)}
\end{center}

\vskip 0.2cm
\begin{center}
{\bf Introduction.}
\end{center}
\vskip 0.2cm

In our article we consider some  algebraical methods which may be
useful in some inverse spectral problems  [1]--[3]. In those papers
problems of  recovery of the coefficients of the boundary conditions
are considered.

In  [1]--[3]  the boundary conditions of the inverse  spectral
problem are determined by the $2\times 4$-matrix $A$ (${\rm rank}\,
A = 2$):

$$ A = \left[ \begin{array}{cccc}
a_{11} & a_{12}          & a_{13} & a_{14} \\
a_{21}        & a_{22} & a_{23} & a_{24}
 \end{array} \right]
$$

 The boundary conditions  which are determined by the $2\times 4$-matrix $A$ and
the boundary conditions  which are determined by the $2\times
4$-matrix $B$ (where  ${\rm  rank}\, A = {\rm rank}\, B = 2$) are
equivalent if and only if their matrices  $A$ and $B$ are linearly
equivalent, i.e. there exists non-singular ${2 \times 2}$-matrix $S$
such that $ A = SB$ (see [3]). Hence, the search for the boundary
conditions  is equivalent to finding the class of linearly
equivalent $2\times 4$-matrices.

Let
$$ A_{ij} = \left| \begin{array}{cc}
a_{1i} & a_{1j}          \\
a_{2i}        & a_{2j}
 \end{array} \right|  (i=1,2,3,4;j=1,2,3,4 )
$$
 are minors of the matrix $A$.

Let's remind $A_{ij}=-A_{ji}$.

By methods from [1]--[3] minors
$A_{12},A_{13},A_{14},A_{23},A_{24},A_{34}$ can be found apart from
a constant.

We want to find a matrix which has minors $A_{12}$, $A_{13}$,
$A_{14}$, $A_{23}$, $A_{24}$, $A_{34}$. More precisely, we want to
find coefficients of the matrix.

\vskip 0.2cm
\begin{center}
\S~1 {\bf Reconstruction of a matrix by its minors.}
\end{center}
\vskip 0.2cm

\vskip 0.2cm

{\bf Theorem~1.} {\it Let $A$, $B$ -- $2\times 4$-matrices
and ${\rm rank}\, A = {\rm rank}\, B = 2$. \\
$
 A_{ij} = \left| \begin{array}{cc}
a_{1i} & a_{1j}          \\
a_{2i}        & a_{2j}
 \end{array} \right|  (i=1,2,3,4; j=1,2,3,4 )
$
are minors of the matrix $A$.\\
$
 B_{ij} = \left| \begin{array}{cc}
b_{1i} & b_{1j}          \\
b_{2i}        & b_{2j}
 \end{array} \right|  (i=1,2,3,4; j=1,2,3,4 )
$ are minors of the matrix $B$.

$\left< {\bf a}_1,{\bf a}_2 \right>$ is a span of the row-vectors
${\bf a}_1= (a_{11}, a_{12}, a_{13}, a_{14})$  and ${\bf a}_2=
(a_{21}, a_{22}, a_{23}, a_{24})$.

$\left< {\bf b}_1,{\bf b}_2 \right>$ is a span of the row-vectors
${\bf b}_1= (b_{11},b_{12}, b_{13}, b_{14})$  and ${\bf b}_2=
(b_{21},b_{22}, b_{23}, b_{24})$.

Then the following statements are equivalent:

i. Exists a  number $t \ne 0$ such that $ A_{ij } = t \, B_{ij}, (i
= 1, 2,3,4; j = 1, 2, 3, 4)$   where $t$ does not depend on the
suffixes $  i, j$ chosen.

ii. $\left< {\bf a}_1,{\bf a}_2 \right> = \left< {\bf b}_1,{\bf b}_2
\right>.$

iii. The matrices  $A$ and $B$ are  linearly equivalent, i.e. there
exist non-singular ${2 \times 2}$-matrix $S$ such that $ B  = S \,
A$
}

\vskip 0.2cm {\bf Proof.}

{\bf i $\Rightarrow$  ii.} The condition that any vector
$(x_1,x_2,x_3,x_4)$ should lie in the linear envelope $\left< {\bf
a}_1,{\bf a}_2 \right>$ is that the matrix
$$
\left[
\begin{array}{cccc}
a_{11} & a_{12}          & a_{13} & a_{14} \\
a_{21}        & a_{22} & a_{23} & a_{24} \\
x_1 & x_2 & x_3 & x_4
 \end{array} \right]
 $$
 should be of $rank \, 2$. In this case the determinant of the
$3\times 3$-submatrix  must be zero. Expanding these determinants in
terms of the last row, we get

$$
 \left\{
\begin{array}{ccccccccc}
A_{23}x_1&+&A_{31}x_2&+&A_{12}x_3& & & = & 0 \\
A_{24}x_1&+&A_{41}x_2& & &+&A_{12}x_4 & = & 0 \\
A_{34}x_1&+& & & A_{41}x_3&+&A_{13}x_4& = & 0 \\
& & A_{34}x_2&+&A_{42}x_3&+&A_{23}x_4& = & 0
 \end{array}  \right. \eqno (1)
$$
Similarly, $(x_1,x_2,x_3,x_4)$ should lie in the linear envelope
$\left< {\bf b}_1,{\bf b}_2 \right>$ if and only if
$$
 \left\{ \begin{array}{ccccccccc}
B_{23}x_1&+&B_{31}x_2&+&B_{12}x_3& & & = & 0 \\
B_{24}x_1&+&B_{41}x_2& & &+&B_{12}x_4 & = & 0 \\
B_{34}x_1&+& & & B_{41}x_3&+&B_{13}x_4& = & 0 \\
& & B_{34}x_2&+&B_{42}x_3&+&B_{23}x_4& = & 0
 \end{array}  \right.\eqno (2)
$$
$A_{ij } = t \, B_{ij}, (i = 1, 2,3,4; j = 1, 2, 3, 4)$   where
$t\ne  0$ is not depend on the suffixes $ i, j$ chosen. Hence, the
systems of the linear equations (1), (2) are equivalent. Therefore
$\left< {\bf a}_1,{\bf a}_2 \right> = \left< {\bf b}_1,{\bf b}_2
\right>.$

 {\bf ii $\Rightarrow$  iii.} Let ${\bf
L}=\left< {\bf a}_1,{\bf a}_2 \right> = \left< {\bf b}_1,{\bf b}_2
\right>.$ $rank \, A = rank \, B = 2$, hence,  ${\bf a}_1,{\bf a}_2$
is basis in ${\bf L}$ and ${\bf b}_1, {\bf b}_2$ is basis in ${\bf
L}$ too. Therefore there exist a non-singular ${2 \times 2}$-matrix
$S$ of basis transformation such that ${\bf b}_i=s_{i1}\, {\bf a}_1
+ s_{i2}\, {\bf a}_2, (i = 1, 2).$ Hence, $ B  = S \, A$.

 {\bf iii $\Rightarrow$ i.} Suppose that $ B =
S \, A$, where $S$ is a  non-singular ${2 \times 2}$-matrix. We
consider the ${2 \times 2}$-submatrix of $A$ by selecting the
columns with suffixes $i, j$ and corresponding submatrix of $B$. We
have
$$
 \left[ \begin{array}{cc}
b_{1i} & b_{1j}          \\
b_{2i}        & b_{2j}
 \end{array} \right]
=
 \left[ \begin{array}{cc}
s_{11} & s_{12}          \\
s_{21}        & s_{22}
 \end{array} \right]
\,
 \left[ \begin{array}{cc}
a_{1i} & a_{1j}          \\
a_{2i}        & a_{2j}
 \end{array} \right]
$$.
and taking determinants of both the sides, we find that $ A_{ij } =
t \, B_{ij}$ , where $t = \det\, S \ne 0$ is not depend on the
suffixes $ i, j$.

The theorem is proved. \vskip 0.2cm

 If we know minors $ A_{ij }, (i
= 1, 2,3,4; j = 1, 2, 3, 4)$ of  the matrix $A$, by (1) we can find
a matrix  which is  linearly equivalent to $A$ .

\vskip 0.2cm

\vskip 0.2cm
\begin{center}
\S~2 {\bf  Examples.}
\end{center}
\vskip 0.2cm

{\bf Example~1.} Let $A_{12}\ne 0$, then  (1) is  equivalent to
$$
\left\{ \begin{array}{ccccccccc}
A_{23}x_1&+&A_{31}x_2&+&A_{12}x_3& & & = & 0 \\
A_{24}x_1&+&A_{41}x_2& & &+&A_{12}x_4 & = & 0
 \end{array}  \right.
$$
And the matrix $A$ is   linearly equivalent to the following  matrix
$$
  \left[ \begin{array}{cccc}
1 & 0 & -A_{23} & -A_{24} \\
0 & 1 &  A_{13} &  A_{14}
 \end{array} \right].
$$
{\bf Example~2.} Let $A_{13}\ne 0$, then  (1) is  equivalent to
$$
\left\{ \begin{array}{ccccccccc}
A_{23}x_1&+&A_{31}x_2&+&A_{12}x_3& & & = & 0 \\
A_{34}x_1&+& & & A_{41}x_3&+&A_{13}x_4& = & 0

 \end{array}  \right.
$$
And the matrix $A$ is   linearly equivalent to the following  matrix
$$
  \left[ \begin{array}{cccc}
1  & A_{23} & 0 & -A_{34} \\
0  &  A_{12} & 1 & A_{14}
 \end{array} \right].
$$
{\bf Example~3.} Let $A_{14}\ne 0$, then  (1) is  equivalent to
$$
\left\{ \begin{array}{ccccccccc}
A_{24}x_1&+&A_{41}x_2& & &+&A_{12}x_4 & = & 0 \\
A_{34}x_1&+& & & A_{41}x_3&+&A_{13}x_4& = & 0

 \end{array}  \right.
$$
And the matrix $A$ is   linearly equivalent to the following  matrix
$$
  \left[ \begin{array}{cccc}
1  & A_{24} & A_{34}& 0  \\
0  &  A_{12} & A_{13}& 1
 \end{array} \right].
$$
{\bf Example~4.} Let $A_{23}\ne 0$, then  (1) is  equivalent to
$$
\left\{ \begin{array}{ccccccccc}
A_{23}x_1&+&A_{31}x_2&+&A_{12}x_3& & & = & 0 \\
& & A_{34}x_2&+&A_{42}x_3&+&A_{23}x_4& = & 0
 \end{array}  \right.
$$
And the matrix $A$ is   linearly equivalent to the following  matrix
$$
  \left[ \begin{array}{cccc}
A_{13}&1&0 & -A_{34}  \\
- A_{12} &0&1& A_{24}
 \end{array} \right].
$$
{\bf Example~5.} Let $A_{24}\ne 0$, then  (1) is  equivalent to
$$
\left\{ \begin{array}{ccccccccc}
A_{24}x_1&+&A_{41}x_2& & &+&A_{12}x_4 & = & 0 \\
& & A_{34}x_2&+&A_{42}x_3&+&A_{23}x_4& = & 0
 \end{array}  \right.
$$
And the matrix $A$ is   linearly equivalent to the following  matrix
$$
  \left[ \begin{array}{cccc}
A_{14}&1&-A_{34}&0  \\
- A_{12} &0& A_{23}&1
 \end{array} \right].
$$
{\bf Example~6.} Let $A_{34}\ne 0$, then  (1) is  equivalent to
$$
\left\{ \begin{array}{ccccccccc}
A_{34}x_1&+& & & A_{41}x_3&+&A_{13}x_4& = & 0 \\
& & A_{34}x_2&+&A_{42}x_3&+&A_{23}x_4& = & 0
 \end{array}  \right.
$$
And the matrix $A$ is   linearly equivalent to the following  matrix
$$
  \left[ \begin{array}{cccc}
A_{14}&A_{24}&1&0  \\
- A_{12} &- A_{23}&0&1
 \end{array} \right].
$$

\vskip 0.2cm
\begin{center}
\S~3 {\bf  Plucker relation.}
\end{center}
\vskip 0.2cm

\vskip 0.2cm

{\bf Theorem~2.} {\it Let $A_{12}, A_{13}, A_{14}, A_{23}, A_{24},
A_{34}$  are some numbers not all  equal to zero. Then the following
statements are equivalent:

i. There exists a $2\times 4$-matrix $A$  such that $A_{12}, A_{13},
A_{14}, A_{23}, A_{24}, A_{34}$
 are  minors of  $A$ and  $rank  \, A =2$.

ii. The following condition is satisfied
$$
 A_{12} \, A_{34} - A_{13} \, A_{24} + A_{14} \, A_{23} = 0.
\eqno (3)
$$
 }
 This condition is called {\rm Plucker relation} (see [4],~[5]).

\vskip 0.2cm

{\bf Proof.}

 {\bf i $\Rightarrow$  ii.} Since
$$
\left| \begin{array}{cccc}
a_{11} & a_{12}          & a_{13} & a_{14} \\
a_{21}        & a_{22} & a_{23} & a_{24} \\
a_{11} & a_{12}          & a_{13} & a_{14} \\
a_{21}        & a_{22} & a_{23} & a_{24}
 \end{array} \right| =0,
$$
by Laplace's expansion of this determinant  we  get
$$
2 \, A_{12} \, A_{34} - 2 \, A_{13} \, A_{24} + 2 \,  A_{14} \, A_{23} = 0.
$$
Hence, Plucker condition (3)   is satisfied.

{\bf ii $\Rightarrow$  iii.} Since ${\rm rank}\, A = 2$, for at
least one set of suffixes $i,j$, $A_{ij}$ is not equal to zero.
Let's assume, for definiteness, that $A_{12} \ne 0$. Then by Plucker
condition (3) we get
$$
  A_{34}= \frac{ A_{13} \, A_{24} -  A_{14} \, A_{23}}{A_{12}} .
$$
Consider the following $2\times 4$-matrix
$$
  \left[ \begin{array}{cccc}
A_{12} & 0 & -A_{23} & -A_{24} \\
0 & 1 & \frac{ A_{13}}{A_{12}} & \frac {A_{14}} {A_{12}}
 \end{array} \right].
$$
It is easy to check up that
$A_{12}, A_{13}, A_{14}, A_{23}, A_{24}, A_{34}$
 are  minors of  this matrix.

The theorem is proved.

\vskip 0.2cm
\begin{center}
\S~4 {\bf  Approximation by a method of orthogonal projection.}
\end{center}
\vskip 0.2cm

When we measure and we calculate something small errors are
possible. Hence it is possible that the numbers $
\widetilde{A}_{12}$, $  \widetilde{A}_{13} $, $\widetilde{A}_{14}$,
$\widetilde{A}_{23}$, $ \widetilde{A}_{24}$, $ \widetilde{A}_{34}$
which are found by methods from [1]--[3] do not satisfy to Plucker
relation (3). Then these numbers are not minors of a matrix.
Therefore we must find numbers $A_{12}$, $A_{13}$, $A_{14}$,
$A_{23}$, $A_{24}$, $A_{34}$ close to the values $
\widetilde{A}_{12}$, $ \widetilde{A}_{13} $, $\widetilde{A}_{14} $,
$\widetilde{A}_{23}$, $ \widetilde{A}_{24}$, $ \widetilde{A}_{34}$
and satisfy to Plucker relation (3).

We have
$$ A_{12} \, A_{34} - A_{13}
\,  A_{24} +  A_{14} \, A_{23} = 0.
$$

By definition, put

$ \, x_1 = A_{12},
 \, x_2 = A_{34},
 \, x_3 = A_{13},
 \, x_4 = -A_{24},
\, x_5 = A_{14},
 \, x_6 = A_{23}.
$ Using this definition, we get the Plucker relation

$$
 x_1\, x_2 + x_3\, x_4+ x_5\, x_6  = 0,
\eqno (4)
$$

which characterizes  a surface $F$ in the 6-dimensional space.

By definition, put

$ \, y_1 = \widetilde{A}_{12},
 \, y_2 = \widetilde{A}_{34},
 \, y_3 = \widetilde{A}_{13},
 \, y_4 = -\widetilde{A}_{24},
\, y_5 = \widetilde{A}_{14},
 \, y_6 = \widetilde{A}_{23}.
$

By definition, put
$$
\begin{array}{c}
X = (x_1, x_2 ,  x_3, x_4,  x_5, x_6) ,\quad Y = (y_1, y_2, y_3,
y_4, y_5, y_6),
\\
X^*=(x_2, x_1, x_4, x_3,  x_6, x_5),\quad Y^*=(y_2, y_1, y_4, y_3,
y_6, y_5),
\\
(\vec X,\vec Y)=x_1\, y_1 + x_2\, y_2+x_3\, y_3 + x_4\, y_4+x_5\,
y_5 + x_6\, y_6, \end{array}
$$
where $\quad \vec X=\overrightarrow{OX}$,  $\vec
Y=\overrightarrow{OY}$, $O$ is the origin.

Let $X$ be an orthogonal projection of $Y$ on  surface  (4). The
vector $\vec X^*$ is  normal for the surface (4) in the point $X$.
 It is identical to the following equations
$$ \vec Y = \vec X + p\,\vec X^*,\eqno (5)
$$
$$
(\vec X, \vec X^*) = 0, \eqno (6)
$$
 where $p$ is a real number.
Having solved a set of linear equations (5) with the unknowns $x_1,
x_2 ,  x_3, x_4,  x_5, x_6$, we obtain

$$
 \vec X = \frac1{1-p^2}\, (\vec Y - p\, \vec Y^*).
\eqno (7)
$$
 From (7) it is easy to obtain
$$
\vec X^* = \frac1{1-p^2}\, (\vec Y^* - p\, \vec Y). \eqno (8)
$$
Substituting (7) for $\vec X$ and (8) for $\vec X^*$ in (6), we
obtain
$$ (\vec Y - p\, \vec Y^*,\, \vec Y^* - p\, \vec Y) = 0.
$$

Notice that
$$
(\vec Y,\vec Y^*)\ne 0,\quad (\vec Y^*,\vec Y^*)=(\vec Y,\vec Y),
\quad (\vec Y^*,\vec Y)=(\vec Y,\vec Y^*).
$$
Therefore,
$$
p^2 - 2\, p\, \frac {(\vec Y,\vec Y)}{(\vec Y,\vec Y^*)} + 1 = 0.
$$
This quadric equation has two roots
$$ p= \frac
{(\vec Y,\vec Y)\mp\sqrt{(\vec Y,\vec Y)^2 - (\vec Y,\vec
Y^*)^2}}{(\vec Y,\vec Y^*)}.
$$
If $X$ is close to  $Y$, then $|p|<<1$ and thus we have
$$
p= \frac {(\vec Y,\vec Y)-\sqrt{(\vec Y,\vec Y)^2 - (\vec Y,\vec
Y^*)^2}}{(\vec Y,\vec Y ^*)}. \eqno (9)
$$
The vector $\vec X$ can be found by using (7) and (9). The
coordinates  $A_{12}$,  $A_{13}$, $A_{14}$, $A_{23}$, $A_{24}$,
$A_{34}$ of $X$ are minors of a matrix. This matrix can be found by
using (1).

\vskip 0.2cm

\vskip 0.3cm \centerline{REFERENCES} \vskip 0.2cm
\begin{enumerate}
\item[{[1]}] {\it Akhatov~I.~Sh.,  Akhtyamov~A.~M.\/}
Determination of the form of attachment of a rod using the natural
frequencies of its flexural oscillations. ~// J. Appl. Math. Mech.,
vol.~65 (2001), no.~2, pp.~283--290.

\item[{[2]}] {\it Akhtyamov~A.~M.\/}
 On uniqueness  of the reconstruction of the boundary conditions of a spectral problem from its  spectrum
(Russian). ~// Fundam. Prikl. Mat, vol.~6 (2000), no.~4,
pp.~995--1006.
\item[{[3]}] {\it Akhtyamov~A.~M., Mouftakhov~A.~V.\/}
 Identification of boundary conditions using
natural frequencies, ~// Inverse Problems in Egineering, vol.~22,
no.~3, 393--408 (2004).
\item[{[4]}] {\it Hodge~W.~V.~D., Pedoe~D.\/}
Methods of algebraic geometry. Vol.~1. University Press, Cambridge,
1994.
\item[{[5]}] {\it Postnikov~M.~M.\/}
Linear Algebra and Differential Geometry. MIR, Moscow, 1982.
\end{enumerate}
}

\end{document}